\documentclass[reqno]{amsart}

\usepackage{times}
\usepackage{mathrsfs}
\usepackage{amsmath}
\usepackage{amsthm}
\usepackage{amsfonts}

\newtheorem{thm}{Theorem}[section]
\newtheorem{lem}[thm]{Lemma}

\newtheorem{prop}[thm]{Proposition}
\newtheorem{rems}[thm]{Remarks}
\newtheorem{rem}[thm]{Remark}
\newtheorem{example}[thm]{Example}

\DeclareMathAlphabet{\mathpzc}{OT1}{pzc}{m}{it}

\numberwithin{equation}{section}

\newcommand{\Wqb}{W_{q,\mathcal{B}}}
\newcommand{\Bqb}{B_{q,p;\mathcal{B}}}
\newcommand{\R}{\mathbb{R}}

\newcommand{\N}{\mathbb{N}}
\newcommand{\F}{\mathbb{F}}
\newcommand{\A}{\mathbb{A}}
\newcommand{\Ac}{\mathcal{A}}
\newcommand{\E}{\mathbb{E}}
\newcommand{\X}{\mathbb{X}}

\newcommand{\ml}{\mathcal{L}}

\newcommand{\mH}{\mathcal{H}}
\newcommand{\K}{\mathcal{K}}

\newcommand{\Om}{\Omega}

\newcommand{\rd}{\mathrm{d}}

\newcommand{\bqn}{\begin{equation}}
\newcommand{\eqn}{\end{equation}}
\newcommand{\bqnn}{\begin{equation*}}
\newcommand{\eqnn}{\end{equation*}}
\newcommand{\bear}{\begin{eqnarray}} 
\newcommand{\eear}{\end{eqnarray}} 
\newcommand{\bean}{\begin{eqnarray*}} 
\newcommand{\eean}{\end{eqnarray*}} 
\newcommand{\bs}{\begin{split}}
\newcommand{\es}{\end{split}}

\newcommand{\dhr}{\mathrel{\lhook\joinrel\relbar\kern-.8ex\joinrel\lhook\joinrel\rightarrow}}

\title[Global Bifurcation of Positive Equilibria]
{Global Bifurcation of Positive Equilibria in Nonlinear Population Models}

\author[Ch. Walker]{Christoph Walker}

\address{Leibniz Universit\"at Hannover, Institut f\"ur Angewandte Mathematik, Welfengarten 1, D--30167 Hannover, Germany.}
\email{walker@ifam.uni-hannover.de}

\begin{document}

\begin{abstract}
Existence of nontrivial nonnegative equilibrium solutions for age structured population models with nonlinear diffusion is investigated. Introducing a parameter measuring the intensity of the fertility, global bifurcation is shown of a branch of positive equilibrium solutions emanating from the trivial equilibrium. Moreover, for the parameter-independent model we establish existence of positive equilibria by means of a fixed point theorem for conical shells. 
\end{abstract}

\keywords{Population models, age structure, nonlinear diffusion, global bifurcation, maximal regularity.
\\
{\it Mathematics Subject Classifications (2000)}: 35B32, 35K55, 35K90, 92D25.}

\maketitle

\section{Introduction}

\noindent This paper deals with finding nonnegative equilibrium solutions to age and spatially structured population equations. In an abstract setting the problem reads: find a nontrivial function $u:J\rightarrow E_0^+$ satisfying the nonlinear problem
\begin{align}
&\partial_au \, +\, A(u,a)\,u\,+\,\mu(u,a)\, u =0\ ,\qquad a\in J\setminus\{0\}\ ,\label{1}\\ 
&u(0)\, =\, \int_0^{a_m}\beta(u,a)\, u(a)\, \rd a\ ,\label{2}
\end{align}
where $E_0$ is an ordered Banach space with positive cone $E_0^+$ and $J:=[0,a_m)$ with $a_m\in (0,\infty]$ denotes the maximal age. Equations \eqref{1}, \eqref{2} arise naturally when considering  equilibrium (i.e. time-independent) solutions to population models, where the function $u$ represents the density of a population of individuals with age $a\in J$ whose evolution is governed by death and birth processes according to the density dependent death modulus $\mu(u,a)$ and birth modulus $\beta(u,a)$, respectively. The term $A(u,a)$ is (for $u$ and $a$ fixed) a linear unbounded operator $A(u,a):E_1\subset E_0\rightarrow E_0$ defined on some common subspace $E_1$ of $E_0$ and models spatial movement of individuals. The full non-equilibrium equations involve an additional time derivative in \eqref{1}, \eqref{2}. Such age-structured equations have been studied since long ago (see \cite{Webb,WebbSpringer} and the references therein), in particular in situations where spatial movement is neglected (i.e. $A\equiv 0$) or when $A$ does not depend on the density $u$ itself (see \cite{LanglaisBusenberg,LanglaisJMB,MagalThieme,RhandiSchnaubelt,ThiemeDCDS98} and the references therein). Less seems to be known about models involving nonlinear age-dependent diffusion (however, see \cite{BusenbergIannelli3,LanglaisSIAM,LauWalJMPA,WalkerEJAM,WalkerDCDS}).  

To understand the asymptotic behavior of the time evolution of structured populations a precise knowledge about equilibrium solutions (i.e. solutions to \eqref{1}, \eqref{2}) is needed. In the present paper we focus on such solutions for the nonlinear age-dependent case $A=A(u,a)$. Clearly, $u\equiv 0$ is a solution to \eqref{1}, \eqref{2} and thus the aim is to give conditions for existence of nontrivial solutions. Moreover, since $u$ in \eqref{1}, \eqref{2} represents a density, any solution should be nonnegative an thus, in the abstract setting, belong to the positive cone $E_0^+$.

To shorten notation we introduce an operator $\A$ as
\bqn\label{a}
\A(u,a):= A(u,a)+\mu(u,a)\ .
\eqn
Suppose that for $u$ fixed, the map $a\mapsto \A(u,a)$ generates on the Banach space $E_0$ a parabolic evolution operator $\Pi_u(a,\sigma)$, $0\le \sigma\le a<a_m$. Then an easy -- but fundamental -- consequence of properties of evolution operators is that any solution $u$ to \eqref{1}, \eqref{2} must satisfy the relation
\bqn\label{3}
u(a)=\Pi_u(a,0)u(0)\ ,\quad a\in J\ ,\qquad u(0)=Q(u)u(0)\ ,
\eqn
where the linear operator $Q(u)$ on $E_0$ is (for $u$ fixed) given by
$$
Q(u):=\int_0^{a_m}\beta(u,a)\Pi_u(a,0)\,\rd a\ .
$$
Roughly speaking, $Q(u)$ contains information about the spatial distribution of the average number of offspring per individual over the entire lifespan of the individual. If spatial movement is neglected, that is, if $A\equiv 0$ and hence $\Pi_u(a,\sigma)=e^{-\int_{\sigma}^a \mu(u,r)\rd r}$, then $Q(u)$ is simply the {\it net reproduction rate} (see \cite{Webb}), and for any solution $u$ to \eqref{1}, \eqref{2} with $A\equiv 0$ this number $Q(u)$ necessarily equals 1 according to \eqref{3}. If spatial movement is included, then \eqref{3} implies that $u(0)$ is (if nonzero) an eigenvector to the eigenvalue 1 of the operator $Q(u)$.\\

In Section \ref{Sec3} we suggest a bifurcation approach to establish positive solutions emanating from the trivial solution $u\equiv 0$. We introduce a bifurcation parameter $n$, which determines the intensity of the fertility without changing its structure, by setting
\bqn\label{4}
n\, b(u,a):=\beta(u,a)\ ,
\eqn
where $b$ is normalized such that the spectral radius $r(Q_0)$ of the bounded linear operator
$$
Q_0:=\int_0^{a_m} b(0,a) \Pi_0 (a,0)\, \rd a
$$
satisfies $r(Q_0)=1$.
Hence $r(Q(0))=n r(Q_0)=n$ and $n$ thus represents the ``inherent net reproduction rate at low densities'' (technically when $u\equiv 0$). If $Q_0$ is a compact positive operator, then 1 is an eigenvalue of $Q_0$ and there is a distribution for which spatial movement, birth, and death processes balance each other if these processes are governed according to $A(0,\cdot)$, $b(0,\cdot)$, and $\mu(0,\cdot)$.

In \cite{WalkerSIMA} a local bifurcation result was shown and a set of $n$-values around the critical value $n=1$ was provided for which \eqref{1}, \eqref{2} subject to \eqref{4} possess nontrivial positive solutions. The aim of this paper is to extend these results to prove a global bifurcation phenomenon. More precisely, if the operator $A(u,a)$ admits a suitable decomposition $A(u,a)=A_0(a)+A_*(u,a)$ with $A_*$ being of ``lower order'' (see Remark~\ref{KK} and \eqref{32}) we use Rabinowitz's alternative \cite{RBA} to show that there is an unbounded continuum of nontrivial positive solutions $(n,u)$ to \eqref{1}, \eqref{2} subject to \eqref{4}. Furthermore, we characterize the set of $n$-values more detailed in some cases and give an example demonstrating that the assumptions imposed are quite natural. We shall point out that our results and methods were inspired by \cite{Cushing2}, where global bifurcation for population models neglecting spatial structure from the outset were investigated. More results on bifurcation for age structured equations with $A\equiv 0$ can be found in \cite{Cushing1,Cushing2,Cushing3} and, respectively, in \cite{DelgadoEtAl,DelgadoEtAl2} for linear and age-independent $A$. We also refer to \cite{Pruess1,Pruess2,Pruess3,Webb} and the references therein for equilibrium solutions for age structured equations in general.\\

The subsequent section is then devoted to a different approach for establishing solutions to the parameter-independent problem \eqref{1}, \eqref{2} not assuming a particular decomposition of $A$. This approach covers ``fully'' quasilinear problems and is more or less independent of the previous considerations. Note that the form of the solution in \eqref{3} allows one to interpret $(u,B)$ with $B:=u(0)$ as a fixed point of the map $$(u,B)\mapsto \big(\Pi_u(\cdot,0)B,Q(u)B\big)\ .$$ In Section \ref{Sec2} we extend an argument of \cite[Thm.1]{Pruess2} (see also \cite[Thm.4.1]{Webb}) for non-diffusive population equations which is based on a fixed point theorem for conical shells \cite{AmannSIAM}. The use of such a theorem prevents hitting the trivial solution $u\equiv 0$ (together with $B=0$ being obviously a fixed point of the map above). We thus prove existence of nontrivial positive solutions for \eqref{1}, \eqref{2} under fairly general assumptions. Loosely speaking, nontrivial positive solutions exist provided that the spectral radius of $Q(u)$ for small populations with density $u$ satisfies $r(Q(u))=r(Q(0))=1$ and is an eigenvalue with a common eigenvector $B$ for all $Q(u)$, and provided that, in addition, for large populations densities there holds $r(Q(u))\le 1$. Thus relevant equilibrium solutions exist if small populations do not affect the spatial distribution of net reproduction rate and large populations have a spatial net reproduction rate not exceeding 1.\\

We conclude the introduction with some notation being used in the following.
If $E$ and $F$ are Banach spaces we write $\ml(E,F)$ for the set of linear bounded operators from $E$ to $F$, and we put $\ml(E):=\ml(E,E)$. The subset thereof consisting of compact operators is denoted by $\K(E,F)$ and $\K(E)$, respectively. We write $r(A)$ for the spectral radius of an operator $A\in\ml(E)$. For an ordered Banach space $E$ we let $\ml_+(E)$ denote the positive linear operators and $\K_+(E)$ is the set of compact positive linear operators. Next, $\mathcal{L}is(E,F)$ stands for the set of topological isomorphisms $E\rightarrow F$. By $E\stackrel{d}{\hookrightarrow}F$ we mean that $E$ is densely embedded in $F$ and $E\dhr F$ stands for a compact embedding of $E$ in $F$. If $E\stackrel{d}{\hookrightarrow}F$ we let $\mH(E,F)$ denote the set of all negative generators of strongly continuous analytic semigroups on $F$ with domain $E$. Moreover, given $\omega>0$ and $\kappa\ge 1$ we write $A\in \mH(E,F;\kappa,\omega)$ if $A\in\ml(E,F)$ is such that $\omega+A\in \mathcal{L}is(E,F)$ and
$$
\frac{1}{\kappa}\,\le\,\frac{\|(\lambda+A)u\|_{F}}{\vert\lambda\vert \,\| u\|_{F}+\|u\|_{E}}\,\le \, \kappa\ ,\quad Re\, \lambda\ge \omega\ ,\quad u\in E\setminus\{0\}\ .
$$
Note that
$$
\mathcal{H}(E,F)=\bigcup_{ \substack{\kappa\ge 1 \\ \omega>0}} \mathcal{H}(E,F;\kappa,\omega)\ .
$$
We refer to \cite{LQPP} for more details. 

Given open subsets $X\subset E$, $Y\subset F$ and another Banach space $G$ we mean by $f\in C^{0,\rho}(X\times Y,G)$ for $\rho>0$ a continuous map $f:X\times Y\rightarrow G$ such that $f(x,\cdot)$ is $\rho$-H\"older continuous for each $x\in X$. Then
$$
\big[f(x,\cdot)\big]_{\rho,Y}:=\sup_{\substack{y_1, y_2\in Y\\ y_1\not= y_2}} \dfrac{\| f(x,y_1)-f(x,y_2)\|_G}{\| y_1-y_2\|_F^\rho}\ .
$$
We let $C_b(E,F)$ denote the continuous functions from $E$ to $F$ being bounded on bounded sets.\\

\subsection*{General assumptions.} 
Throughout the paper we assume that
$E_0$ is a real Banach space ordered by a closed convex cone $E_0^+$ and $E_1\stackrel{d}{\dhr}E_0$ for some Banach space $E_1$.
We fix for each $\theta\in (0,1)$ an admissible interpolation functor $(\cdot,\cdot)_\theta$, that is, an interpolation functor $(\cdot,\cdot)_\theta$ such that $E_1\stackrel{d}{\hookrightarrow}E_\theta:=(E_0,E_1)_\theta$. Note that $E_\theta\dhr E_\vartheta$ for $0\le\vartheta<\theta\le 1$ (see \cite[I.Thm.2.11.1]{LQPP}). The interpolation spaces $E_\theta$ are given their natural order induced by the cone $E_\theta^+:=E_\theta\cap E_0^+$. We fix $a_m\in (0,\infty]$ and set $J:=[0,a_m)$. Observe that $a_m=\infty$ is explicitly allowed.

\section{Global Bifurcation of Positive Equilibria}\label{Sec3}

We focus our attention on the parameter dependent problem \eqref{1}, \eqref{2} subject to \eqref{4} in a more general framework. More precisely, we look for solutions $(n,u)$ to problems of the form
\bqn
\begin{aligned}\label{30}
 \partial_au \, +\,     \A(u,a)\,u\, &=0\ ,\qquad a\in J\ ,\\ 
u(0)\,& =\,n\,\ell(u)\ .
\end{aligned}
\eqn
The linear unbounded operator $\A(u,a):E_1\subset E_0\rightarrow E_0$ (for $u$ and $a$ fixed) and the operator $\ell$ with $\ell(u)\in E_0$ and $\ell (0)=0$ are supposed to satisfy some technical assumptions specified later on. The main example for the latter we have in mind is, of course,
$$
\ell(u)=\int_0^{a_m} b(u,a)u(a)\,\rd a\ .
$$
Clearly, the branch $(n,u)=(n,0)$, $n\in\R$, consists of (trivial) solutions to \eqref{30}. Our aim is to prove that another unbounded branch of nontrivial positive solutions $(n,u)$ (i.e. $u(a)\in E_0^+$ for $a\in J$ and $u\not\equiv 0$) bifurcates from the trivial branch at some critical value, which we may assume to be $n=1$ under a suitable normalization. Imposing maximal $L_p$-regularity for (a part of) the operator $\A$ we will show that the nontrivial branch $(n,u)$ is unbounded in $\R^+\times (L_p^+(J,E_1)\cap W_p^1(J,E_0))$ for some $p\in (1,\infty)$ fixed. The result is inspired by the work of \cite{Cushing2} for the non-diffusive case $A\equiv 0$ and is a consequence of Rabinowitz's alternative \cite{RBA}. The application of this alternative requires some compactness in appropriate spaces of the operators involved, which is guaranteed, for example, by a generalized Aubin-Dubinskii lemma (see Remark \ref{KK} and Lemma \ref{A6}).\\

To state the precise assumptions let us fix $p\in (1,\infty)$ and introduce the spaces
$$
\E_0:=L_p(J,E_0)\ ,\qquad \E_1:=L_p(J,E_1)\cap W_p^1(J,E_0)\ .
$$
Recall the embedding
\bqn\label{b}
\E_1\hookrightarrow BUC(J,E_\varsigma)\ ,
\eqn
where $BUC$ stands for bounded and uniformly continuous and $E_\varsigma:=(E_0,E_1)_{\varsigma,p}$ with $(\cdot,\cdot)_{\varsigma,p}$ being the real interpolation functor for $\varsigma:=\varsigma(p):=1-1/p$. The trace $\gamma u:=u(0)$ is thus well-defined for $u\in \E_1$.  We further fix Banach spaces $\F_1$, $\F_2$, $\F_3$, and $\F_4$ such that
\bqn\label{K}
\E_1\dhr \F_j\ ,\quad j=1,2,3,4\ ,
\eqn
and first remark the following:

\begin{rem}\label{KK}
Given $\alpha\in [0,1)$ and $s\in [0,1-\alpha)$, Sobolev spaces $W_p^s(J,E_\alpha)$ are  appropriate choices for $\F_j$ to satisfy \eqref{K}. This follows from the compact embedding $E_1\dhr E_\alpha$ and a generalized Aubin-Dubinskii lemma \cite[Thm.1.1]{AmannGlasnik00} (together with a diagonal sequence argument if $J=\R^+$).
\end{rem}

For the nonlinear operator $\A$ in \eqref{30} we then shall assume a decomposition of the form
\bqn\label{30A}
\A(u,a)=\A_0(a)+\A_*(u,a)\ ,
\eqn
where $\A_0$ is an age-dependent parabolic operator and the nonlinearity of $\A$ in $u$ is contained in a ``lower order perturbation'' $\A_*$. 
To be more precise we suppose for the linear part $\A_0$ that
\bqn\label{31}
\begin{aligned}   & \A_0\in L_{\infty}(J,\ml(E_1,E_0))\ \text{generates a positive parabolic evolution operator}\\
&\Pi_0(a,\sigma), 0\le\sigma\le a<a_m,\ \text{on}\ E_0\ \text{with regularity subspace}\ E_1\ \text{and}\\
&\text{possesses maximal}\ L_p\text{-regularity, that is,}\ (\partial_a+\A_0,\gamma)\in \mathcal{L}is(\E_1,\E_0\times E_{\varsigma})\ ,
\end{aligned}
\eqn
while for the nonlinear part $\A_*$ we assume that
\bqn\label{32}
\A_*\in C\big(\F_1,\ml(\F_2,\E_0)\big)\ \text{with}\ \A_*(0,\cdot)=0\ .
\eqn
We also assume $\ell$ in \eqref{30} admits a decomposition
\bqn\label{33}
\ell(u)=\ell_0(u)+\ell_*(u)
\eqn
with linear part 
\bqn\label{34}
\ell_0\in\ml(\E_1,E_\vartheta)+\ml(\F_3,E_\varsigma) \ \ \text{for some}\ \ \vartheta\in (\varsigma,1]\ ,
\eqn
and nonlinear part 
\bqn\label{35} 
\ell_*\in C_b(\E_1,E_\vartheta)+ C(\F_4,E_\varsigma)\ \ \text{with}\ \ \|\ell_*(u)\|_{E_\varsigma}=o(\|u\|_{\E_1})\ \ \text{as}\ \ \|u\|_{\E_1}\rightarrow 0\ .
\eqn
A consequence of \eqref{31} is that for each datum $(u^0,f)\in E_\varsigma
\times \E_0$ the problem
$$
\partial_a u+\A_0(a)u=f(a)\ ,\quad a\in J\ ,\qquad u(0)=u^0
$$
possesses a unique solution $u\in \E_1$ given by
$$
u(a)=\Pi_0(a,0)u^0+ (K_0f)(a)\ ,\quad a\in J\ ,
$$
satisfying for some number $c_0>0$ independent of $f$ and $u^0$
$$
\| u\|_{\E_1}\le c_0\big(\| u^0\|_{E_\varsigma}+\|f\|_{\E_0}\big)\ ,
$$
where
$$
(K_0f)(a):=\int_0^a\Pi_0(a,\sigma) f(\sigma)\,\rd \sigma\ ,\qquad a\in J\ ,\quad f\in\E_0\ .
$$
Therefore,
\bqn\label{36}
\Pi_0(\cdot,0)\in\ml(E_\varsigma,\E_1)\ ,\qquad K_0\in\ml(\E_0,\E_1)\ ,
\eqn
and \eqref{K}, \eqref{34} thus imply
\bqn\label{uu}
[f\mapsto \ell_0(K_0f)]\in\ml(\E_0,E_\varsigma)\ .
\eqn
We also note that
\bqn\label{38}
Q_0\in \K(E_\varsigma)\quad\text{for}\quad Q_0w:=\ell_0(\Pi_0(\cdot,0)w)\ ,\quad w\in E_\varsigma\  ,
\eqn
which is a consequence of \eqref{36} and either $E_\vartheta\dhr E_\varsigma$ if $\ell_0\in\ml(\E_1,E_\vartheta)$ or \eqref{K} if $\ell_0\in\ml(\F_3,E_\varsigma)$.
Without loss of generality we may assume that $\ell_0$ is normalized such that the spectral radius of $Q_0\in\K(E_\varsigma)$ equals 1, that is, $r(Q_0)=1$.
We first consider the linearization of \eqref{30} around $u\equiv 0$.

\begin{lem}\label{A5}
Suppose \eqref{K}-\eqref{35} and let $r(Q_0)=1$. Then, for each $(h_1,h_2)\in E_\varsigma\times \E_0$, the problem
\bqn\label{776}
\partial_au \, +\,  \A_0(a)\,u\, =h_2(a)\ ,\quad a\in J\ ,\qquad
u(0)-\frac{1}{2}\ell_0(u)\, =\,h_1
\eqn
admits a unique solution $u=S(h_1,h_2)\in\E_1$ given by $
u=\Pi_0(\cdot,0)w(h_1,h_2)+K_0h_2$,
where 
\bqn\label{778}
w(h_1,h_2):=\left(1-\dfrac{1}{2}Q_0\right)^{-1}\left(\dfrac{1}{2}\ell_0(K_0h_2)+h_1\right)\in E_\varsigma\ .
\eqn
The solution operator $S$ belongs to $\ml(E_\varsigma\times\E_0,\E_1)$.
\end{lem}

\begin{proof}
By the previous observations problem \eqref{776} is, for any given $(h_1,h_2)\in E_\varsigma\times \E_0$, equivalent to
\bqnn
\begin{aligned}
&u(a)=\Pi_0(a,0)u(0)+(K_0h_2)(a)\ ,\quad a\in J\ ,\\
&u(0)=\dfrac{1}{2}Q_0u(0)+\frac{1}{2}\ell_0(K_0h_2)+h_1\ .
\end{aligned}
\eqnn
Taking $2>r(Q_0)$ into account, the latter equality entails \eqref{778}, and
defining
$$
S(h_1,h_2):=\Pi_0(\cdot,0)w(h_1,h_2)+K_0h_2\ ,
$$
we derive $S\in\ml(E_\varsigma\times\E_0,\E_1)$ from \eqref{36} and \eqref{uu}.
\end{proof}

Lemma \ref{A5} allows a reformulation of problem \eqref{30} in terms of the operator $S$: writing $n=\lambda+1/2$, a function $u\in\E_1$ solves \eqref{30} if and only if
\bqn\label{400}
u=\lambda Lu+H(\lambda,u)\ ,
\eqn
where
\begin{align*}
Lu:=S\big(\ell_0(u),0\big)\ , \quad H(\lambda,u):=S\big((\lambda+1/2)\ell_*(u),-\A_*(u,\cdot)u\big)\ .
\end{align*}

The maps $L$ and $H$ enjoy the following properties.

\begin{lem}\label{A6}
Suppose \eqref{K}-\eqref{35} and let $r(Q_0)=1$. Then $L\in\K(\E_1)$. Moreover, $H\in C(\R\times\E_1,\E_1)$ is compact and $\|H(\lambda,u)\|_{\E_1}=o(\|u\|_{\E_1})$ as $\|u\|_{\E_1}\rightarrow 0$ uniformly on bounded $\lambda$ intervals.
\end{lem}

\begin{proof}
Clearly, Lemma \ref{A5}, \eqref{K}, and \eqref{34} ensure $L\in\ml(\E_1)$. Using $E_\vartheta\dhr E_\varsigma$ if $\ell_0\in\ml(\E_1,E_\vartheta)$ or $\E_1\dhr \F_3$ if $\ell_0\in\ml(\F_3,E_\varsigma)$, the compactness of $L$ is obvious. Next, Lemma \ref{A5} together with \eqref{32} and \eqref{35} imply $H\in C(\R\times\E_1,\E_1)$. As for its compactness we note that if $(u_j)$ is a bounded sequence in $\E_1$, then $(\ell_*(u_j))$ is relatively compact in $E_\varsigma$ either because $(\ell_*(u_j))$ is bounded in $E_\vartheta\dhr E_\varsigma$  if $\ell_*\in C_b(\E_1,E_\vartheta)$ or because $\E_1\dhr \F_4$ if $\ell_*\in C(\F_4,E_\varsigma)$. Next, \eqref{K} and \eqref{32} ensure that $(\A_*(u_j)u_j)$ is relatively compact in $\E_0$. Hence the compactness of $H$ follows from Lemma \ref{A5}. Finally, the last assertion is a consequence of \eqref{32}, \eqref{35}, and again Lemma \ref{A5}.
\end{proof}

\begin{rem}\label{r33}
Alternatively to \eqref{32} we could have assumed $\A_*\in C(\E_1,\ml(\E_1,\E_0))$ with \mbox{$\A_*(0,\cdot)=0$} is such that $(u\mapsto \A_*(u)u)\in C(\E_1,\E_0)$ is compact.
\end{rem}

To problem \eqref{400} we may now apply Rabinowitz's alternative \cite{RBA}. Recall that the {\it characteristic values} of a linear operator are the reciprocals of its real nonzero eigenvalues. A {\it continuum} in $\R\times\E_1$ is a closed connected subset thereof, and it {\it meets infinity} if it is unbounded.

\begin{prop}\label{A7}
Suppose \eqref{K}-\eqref{35} and let $r(Q_0)=1$ be a simple eigenvalue of $Q_0\in\K(E_\varsigma)$ with eigenvector $B\in E_\varsigma$. Then there exists a maximal continuum $\mathfrak{C}$ in $\R\times\E_1$ consisting of solutions $(n,u)$ to \eqref{30} with $u\not\equiv 0$ if $n\not=1$ and $(1,0)\in\mathfrak{C}$. For $(n,u)\in\mathfrak{C}$ near $(1,0)$ we have
\bqn \label{50}
u=\varepsilon\Pi_0(\cdot,0)B+u_*(\varepsilon)\quad\text{with}\quad \|u_*(\varepsilon)\|_{\E_1}=o(\varepsilon)\quad\text{as}\quad \varepsilon\rightarrow 0\ .
\eqn
The continuum $\mathfrak{C}$ satisfies the following alternative: either $\mathfrak{C}$ meets infinity or it meets a point $(\hat{n},0)$ with \mbox{$\hat{\mu}=\hat{n}-1/2$} being a characteristic value of $L$.
\end{prop}

\begin{proof}
First note that $u=\mu Lu$ with $u\in\E_1$ and $\mu\in\R$ is equivalent to
\bqn\label{43}
u=\Pi_0(\cdot,0)u(0)\ ,\quad u(0)=(\mu+1/2)Q_0 u(0)
\eqn
owing to Lemma \ref{A5} and \eqref{b}. Thus, $\mu+1/2$ is a characteristic value of $Q_0\in \K (E_\varsigma)$ if and only if $\mu$ is a characteristic value of $L\in\K (\E_1)$. Additionally assuming $\mu+1/2$ to be a simple characteristic value of $Q_0$ we claim that $\mu$ is a simple characteristic value of $L$. For, let $u\in\mathrm{ker}\big((\mu L-1)^2\big)\subset \E_1$ and define \mbox{$v:=(\mu L-1)u\in \mathrm{ker}(\mu L-1)\subset\E_1$}. Then \eqref{43} ensures $v=\Pi_0(\cdot,0)v(0)$ with $v(0)$ belonging to $\mathrm{ker}(1-(\mu+1/2)Q_0)$. The characteristic value $\mu+1/2$ of $Q_0$ being simple, we deduce $v(0)=r\xi_0$ for some $r\in\R$ and $\xi_0\in E_\varsigma$ with $\mathrm{ker}(1-(\mu+1/2)Q_0)=\mathrm{span}\{\xi_0\}$. Hence
\bqn\label{44}
v=r\Pi_0(\cdot,0)\xi_0\ ,
\eqn
and we aim for $r=0$. Clearly, from \eqref{43} we have $\Pi_0(\cdot,0)\xi_0\in\mathrm{ker}(\mu L-1)$, whence from \eqref{44} $u=L(\mu u-r\mu\Pi_0(\cdot,0)\xi_0)$. Lemma \ref{A5} then entails
$$
\big(1-(\mu+1/2)Q_0\big) u(0)= -r\mu Q_0\xi_0=\frac{-r\mu}{\mu+1/2}\xi_0\ .
$$
Consequently,
$$
r\mu\xi_0\in \mathrm{rg}\big(1-(\mu+1/2)Q_0\big) \cap\mathrm{ker}\big(1-(\mu+1/2)Q_0\big)=\{0\}
$$
since $\mu+1/2$ is a simple characteristic value of the compact operator $Q_0$. We conclude $r=0$ as desired because $\mu=0$ is impossible owing to the fact that $1/2$ is no characteristic value of $Q_0$ since $r(Q_0)=1$. But then $v=0$ by \eqref{44} and so 
$\mathrm{ker}\big((\mu L-1)^2\big)\subset \mathrm{ker}(\mu L-1)$. Therefore, $\mu$ is indeed simple for $L$ provided $\mu+1/2$ is simple for $Q_0$. In particular, $\mu=1/2$ is a simple characteristic value of $L$ due to the assumption on $r(Q_0)$.

Taking Lemma \ref{A6} into account we may apply Theorem 1.3 from \cite{RBA} to conclude the existence of a maximal continuum $\mathfrak{C}$ in $\R\times \E_1$ with $(1,0)\in\mathfrak{C}$ such that $(n,u)\in\mathfrak{C}$ solves $u=\lambda Lu+ H(\lambda,u)$ with $\lambda=n-1/2$, where $u\not\equiv 0$ if $n\not=1$, and $\mathfrak{C}$ either meets infinity or meets a point $(\hat{n},0)$ with a characteristic value $\hat{\mu}=\hat{n}-1/2$ of $L$ different from $1/2$. Finally, \cite[Lem.1.24]{RBA} implies \eqref{50} and the statement follows.
\end{proof}

\noindent Imposing further conditions on $\A$ and $\ell$ we now prove that a global branch of positive solutions $(n,u)$ to \eqref{30} exists emanating from the trivial branch $(n,0)$, $n\in\R$ at the critical point $(1,0)$. 
To prove this result we suppose that
\bqn\label{51}
\begin{aligned}
&\text{for each}\ u\in\E_1, \A(u,\cdot)\ \text{generates a positive parabolic evolution operator}\\
&\Pi_u(a,\sigma), 0\le\sigma\le a<a_m ,\ \text{on}\ E_0\ \text{with regularity subspace}\ E_1\ .
\end{aligned}
\eqn
We also assume that
\bqn\label{53}
\begin{aligned}
&\text{there is}\ \bar{\ell}_*\ \text{with}\ \ell_*(u)=\bar{\ell}_*(u,u)\ \text{and}\ \bar{\ell}_*(0,\cdot)=0\ \text{such}\\
&\text{that}\ Q_u\in\K_+(E_\varsigma) \text{ for each} \ u\in\E_1\ ,\ \text{where}\\ &Q_uw:=\ell_0(\Pi_u(\cdot,0)w)+\bar{\ell}_*(u,\Pi_u(\cdot,0)w) , w\in E_\varsigma\ .\\
\end{aligned}
\eqn
Note that this definition of $Q_u$ is consistent with \eqref{38}. Let $\mathrm{int}(E_\varsigma^+)$ denote the interior of the positive cone $E_\varsigma^+$. Then we assume further that
\bqn\label{54}
\begin{aligned}
&\text{for each}\ u\in \E_1,\ \text{any positive eigenvector to a}\ \text{positive eigenvalue of}\ Q_u\ \text{belongs to}\ \mathrm{int}(E_\varsigma^+)\ .
\end{aligned}
\eqn
This last assumption is crucial for positivity of solutions but not too restrictive in applications as noted in the following remark (see also Example \ref{ex}).

\begin{rem}\label{A8}
If $\mathrm{int}(E_\varsigma^+)\not=\emptyset$ and $Q_u\in\K_+(E_\varsigma)$ is irreducible (e.g. if strongly positive), then the Krein-Rutman theorem \cite[Thm.12.3]{DanersKochMedina} ensures that the spectral radius $r(Q_u)>0$ is a simple eigenvalue of $Q_u$ with an eigenvector belonging to $\mathrm{int}(E_\varsigma^+)$, and it is the only eigenvalue with positive eigenvector. Thus \eqref{54} holds in this case.
\end{rem}

\begin{thm}\label{A9}
Suppose \eqref{K}-\eqref{35}, and \eqref{51}-\eqref{54}. Further let $r(Q_0)=1$ be a simple eigenvalue of $Q_0\in\K_+(E_\varsigma)$ with eigenvector $B\in \mathrm{int}(E_\varsigma^+)$ and suppose $Q_0$ has no other eigenvalue with positive eigenvector. Then there is a contiunuum $\mathfrak{C}^+$ in $\R^+\times\E_1^+$ of positive solutions to \eqref{30} connecting $(1,0)$ with infinity. 
\end{thm}

\begin{proof}
Let $\mathfrak{C}$ denote the maximal continuum of solutions to \eqref{30} provided by Proposition \ref{A7}. Clearly, if $(n,u)\in \mathfrak{C}$, then $u=\Pi_u(\cdot,0)u(0)$ by \eqref{30}, \eqref{400}, and \eqref{51}. Thus $u\in \E_1^+$ provided $u(0)\in E_\varsigma^+$. Also note that $(0,u)\in \mathfrak{C}$ would imply $u\equiv 0$ and is thus impossible. Due to $B\in \mathrm{int}(E_\varsigma^+)$ it follows from \eqref{b} and \eqref{50} that for $(n,u)\in \mathfrak{C}$ near $(1,0)$ we have for sufficiently small $\varepsilon>0$ 
$$
\frac{1}{\varepsilon}u(0)=B+\frac{1}{\varepsilon}\gamma u_*(\varepsilon)\in E_\varsigma^+\ .
$$
Since $\mu=1/2$ is a simple characteristic value of $L$, we may refer to \cite[Thm.1.40]{RBA} to deduce that $\mathfrak{C}$ is the union of two subcontinua $\mathfrak{C}^{\pm}$, where $\mathfrak{C}^+$ consists of positive solutions near $(1,0)$, and $\mathfrak{C}^+$ (and also $\mathfrak{C}^-$) meets $(1,0)$ and either meets infinity in $\R\times\E_1$ or a point $(\hat{\mu}+1/2,0)$ with $\hat{\mu}$ being a characteristic value of $L$ different from $1/2$. Consequently, $\mathfrak{C}^+\cap (\R^+\times\E_1^+)\not=\emptyset$, and we now show that $\mathfrak{C}^+$ leaves $\R^+\times\E_1^+$ only at the bifurcation point $(n,u)=(1,0)$. For, suppose the continuum $\mathfrak{C}^+$ leaves $\R^+\times\E_1^+$ at some point $(n_*,u_*)$. Then there are $(n_j,u_j)\in \mathfrak{C}^+\cap (\R^+\times\E_1^+)$ with
\bqn\label{61}
(n_j,u_j)\rightarrow (n_*,u_*)\quad\text{in}\quad \R\times\E_1\ .
\eqn
In particular, writing $n_j=\lambda_j+1/2$ we have
$u_j=\lambda_j L u_j+H(\lambda_j,u_j)$, $j\in\N$,
and letting $j\rightarrow\infty$ we obtain from Lemma \ref{A6} that
$
u_*=\lambda_* L u_*+H(\lambda_*,u_*)
$ for $\lambda_*:=n_*-1/2$. Hence
$$\partial_a u_*+\A(u_*,a)u_*=0\ ,\quad u_*(0)=n_*\ell(u_*)
$$
from which
$$u_*=\Pi_{u_*}(\cdot,0)u_*(0)\ ,\quad u_*(0)=n_* Q_{u_*}u_*(0)
$$
by \eqref{51}. Therefore, either $u_*(0)=0$ and then $u_*\equiv 0$ or $u_*(0)\in E_\varsigma^+\setminus\{0\}$ by \eqref{b} in which case $n_*>0$ must be a characteristic value of $Q_{u_*}$ with a positive eigenvector. Thanks to \eqref{54} we derive
\bqn\label{62}
u_*\equiv 0\qquad\text{or}\qquad u_*(0)\in\mathrm{int}(E_\varsigma^+)\ .
\eqn
First suppose that $u_*\equiv 0$. Then $(n_j,u_j)\rightarrow (n_*,0)$ in $\R^+\times\E_1^+$, and we claim that $\lambda_*=n_*-1/2$ is a characteristic value of $L$. 
Indeed, putting $v_j:=\|u_j\|_{\E_1}^{-1}u_j$ and taking into account \eqref{400}, the compactness of the operator $L$, and the fact that $\|H(\lambda,u)\|_{\E_1}=o(\|u\|_{\E_1})$ as $\|u\|_{\E_1}\rightarrow 0$ we may extract a subsequence of $(v_j)$ converging in $\E_1$ toward some $v\in\E_1^+\setminus\{0\}$ with $v=\lambda_* Lv$. Thus $\lambda_*$ is indeed a characteristic value of $L$. As in \eqref{43} this implies that $n_*=\lambda_*+1/2$ is a characteristic value of $Q_0$ with a positive eigenvector $v(0)\in E_\varsigma^+$, whence $n_*=1$ by assumption. Therefore, $\mathfrak{C}^+$ leaves $\R^+\times\E_1^+$ at the bifurcation point $(n_*,u_*)=(1,0)$. Now suppose that $u_*(0)\in\mathrm{int}(E_\varsigma^+)$. Since $\mathfrak{C}^+$ is connected and leaves $\R^+\times\E_1^+$ at $(n_*,u_*)$, there is a sequence $(\bar{n}_j,\bar{u}_j)\in  \mathfrak{C}^+$ with $\bar{u}_j\not\in\E_1^+$ and $(\bar{n}_j,\bar{u}_j)\rightarrow (n_*,u_*)$ in $\R\times\E_1$. According to \eqref{b} we find $m\in\N$ with $\bar{u}_m(0)\in E_\varsigma^+$. But $(\bar{n}_m,\bar{u}_m)\in \mathfrak{C}^+$ and thus
$$
\partial_a\bar{u}_m+\A(\bar{u}_m,a)\bar{u}_m=0\ ,\quad \bar{u}_m(0)=\bar{n}_m\ell(\bar{u}_m)\ ,
$$
that is, $\bar{u}_m=\Pi_{\bar{u}_m}(\cdot,0)\bar{u}_m(0)\in\E_1^+$ due to \eqref{51} contradicting the choice of the sequence $(\bar{u}_j)$. Therefore, $u_*(0)\in\mathrm{int}(E_\varsigma^+)$ is impossible. We have thus shown that $\mathfrak{C}^+$ leaves $\R^+\times\E_1^+$ only at the bifurcation point $(n_*,u_*)=(1,0)$.

It remains to prove that $\mathfrak{C}^+$ does not meet a point $(\hat{\mu}+1/2,0)$ with $\hat{\mu}\not=1/2$ being a characteristic value of $L$. For, suppose 
$\mathfrak{C}^+$  meets such a point. Then we find $(n_j,u_j)\in \mathfrak{C}^+\subset \R^+\times\E_1^+$ with \mbox{$(n_j,u_j)\rightarrow (\hat{\mu}+1/2,0)$} in $ \R^+\times\E_1^+$. Exactly as above one shows that then $\hat{\mu}=1/2$ which is ruled out by assumption. This proves the theorem.
\end{proof}

If the interior of the positive cone $E_\varsigma^+$ is empty, one can prove another result if we put some symmetry conditions on $\A$ and $\ell$. For the eigenvector $B$ of $Q_0$ we only assume that $B\in E_0^+$.

\begin{thm}\label{A11}
Suppose \eqref{K}-\eqref{35}, \eqref{51}, and \eqref{53}. Moreover, suppose that for each $u\in\E_1$, any positive eigenvalue of $Q_u$ has geometric multiplicity 1 and possesses a positive eigenvector $B_u\in E_0^+$. In addition, let $\A_*(-u,\cdot)=\A_*(u,\cdot)$ and $\bar{\ell}_*(-u,\cdot)=\bar{\ell}_*(u,\cdot)$ for $u\in\E_1$. If $r(Q_0)=1$ is a simple eigenvalue of $Q_0\in\K_+(E_\varsigma)$ with an eigenvector $B\in E_0^+$, and if there is no other eigenvalue of $Q_0$ with positive eigenvector, then
$$
\mathfrak{C}_+:=\big\{(n,u)\in \mathfrak{C}\,;\, u(0)\in E_0^+\big\}\cup \big\{(n,-u)\,;\, (n,u)\in  \mathfrak{C}\,,\, u(0)\not\in E_0^+\,,\, n>0\big\}
$$
is an unbounded closed subset of $\R\times\E_1^+$ such that $
\mathfrak{C}_+\setminus\{(1,0)\}$ consists of positive nontrivial solutions to \eqref{30}, where $\mathfrak{C}$ is the maximal continuum from Proposition \ref{A7}.
\end{thm}

\begin{proof}
As in the proof of Theorem \ref{A9} we have $u\in\E_1^+$ for $(n,u)\in \mathfrak{C}$ provided $u(0)\in E_0^+$. If $(n,u)\in \mathfrak{C}$ is such that $u(0)\not\in E_0^+$ and $n>0$, then $u=\Pi_u(\cdot,0)u(0)$ and $u(0)=nQ_u u(0)$. Hence $n>0$ is a characteristic value of $Q_u$ and by assumption there is a corresponding eigenvector $B_u\in E_0^+$ such that $u(0)=r_u B_u$ for some $r_u<0$. Due to the symmetry conditions put on $\A_*$ and $\ell_*$ it is easily seen that $v:=-u$ satisfies $v=\Pi_v(\cdot,0)v(0)$ and $v(0)=n Q_v v(0)$ with $v(0)=-r_uB_u\in E_0^+$. Consequently, $(n,v)\in \R^+\times\E_1^+$ solves \eqref{30}. Hence $\mathfrak{C}_+$ consists of nonnegative solutions only. To show that $\mathfrak{C}$ is unbounded assume to the contrary that $\mathfrak{C}$ meets a point $(\hat{\mu}+1/2,0)$ with a characteristic value $\hat{\mu}$ of $L$ different from $1/2$. Let $(\lambda_j+1/2,u_j)\in \mathfrak{C}$ be such that $(\lambda_j,u_j)\rightarrow (\hat{\mu},0)$ in $\R\times \E_1$. Recall \eqref{b} and set $\bar{u}_j:=u_j$ if $u_j(0)\in E_0^+$ and $\bar{u}_j:=-u_j$ if $u_j(0)\not\in E_0^+$. Then $\bar{u}_j\in \E_1^+$ and $(\lambda_j,\bar{u}_j)\rightarrow (\hat{\mu},0)$ in $\R\times\E_1$. Exactly as in the proof of Theorem \ref{A9} we deduce that $\hat{\mu}+1/2$ is a characteristic value of $Q_0$ with an eigenvector in $E_\varsigma^+$, whence $\hat{\mu}+1/2=r(Q_0)=1$ by assumption contradicting $\hat{\mu}\not=1/2$. Therefore, $\mathfrak{C}$ is unbounded according to Proposition \ref{A7} and so is $\mathfrak{C}_+$. That the latter is closed in in $\R\times\E_1$ is a consequence of the continuity of $L$ and $H$. This concludes the proof.
\end{proof}

Let us consider the set of possible parameter values in more detail. For, suppose the conditions of Theorem \ref{A9}, let $\mathfrak{C}^+$ denote the unbounded continuum in $\R^+\times\E_1^+$ consisting of positive solution to \eqref{30}, and recall that $Q_u\in\K_+(E_\varsigma)$ for $u\in\E_1$ was defined by
$$
Q_uw:=\ell_0(\Pi_u(\cdot,0)w)+\bar{\ell}_*(u,\Pi_u(\cdot,0)w)\ ,\quad w\in E_\varsigma\ ,\quad u\in \E_1\ .
$$
Observe that for $(n,u)\in \mathfrak{C}^+$ we have $u=\Pi_u(\cdot,0)u(0)$ with $u(0)= n Q_u u(0)$,
whence $nr(Q_u)\ge 1$ for $(n,u)\in \mathfrak{C}^+$. In addition, if
\bqn\label{71}
\begin{aligned}
\text{for each}\ u\in \E_1^+, Q_u\in\K_+(\E_\varsigma)\ \text{has only}\ r(Q_u)>0\ \text{as eigenvalue with positive eigenvector}\ ,
\end{aligned}
\eqn
which holds e.g. if the Krein-Rutman theorem applies to $Q_u$, then necessarily 
\bqn\label{72}
nr(Q_u)=1\ ,\quad (n,u)\in \mathfrak{C}^+\ .
\eqn
This observations guarantees a more precise characterization of the spectrum
$$
\sigma:=\big\{n\,;\, \text{there is}\ u\in\E_1^+\ \text{with}\ (n,u)\in\mathfrak{C}^+\setminus\{(1,0)\}\big\}
$$
as well as of the solution set
$$
\Gamma:=\big\{u\,;\, \text{there is}\ n\in\R\ \text{with}\ (n,u)\in\mathfrak{C}^+\setminus\{(1,0)\}\big\}\ .
$$

The next proposition is in the spirit of \cite{Cushing2} for the spatially homogeneous case $A\equiv 0$ in \eqref{a} and the easy proofs carry over to the present situation almost verbatim. We nevertheless include them here for the reader's ease.

\begin{prop}\label{A20}
Suppose the conditions of Theorem \ref{A9} and \eqref{71}. Setting
$$
\sigma_i:=\inf \sigma\ ,\quad \sigma_s:=\sup \sigma
\ ,\quad N_i:=\inf_{u\in\Gamma} r(Q_u)\ ,\quad N_s:=\sup_{u\in\Gamma} r(Q_u)\ ,
$$
we have:
\begin{itemize}
\item[(i)] $0\not\in \sigma\subset\R^+$ and $\Gamma\subset\E_1^+$.\\
\item[(ii)] $0\le\sigma_i\le 1\le\sigma_s\le\infty$.\\
\item[(iii)] $N_i=0$ iff $\sigma_s=\infty$, and if $N_i>0$, then $\sigma_s=1/N_i<\infty$.\\
\item[(iv)] $N_s=\infty$ iff $\sigma_i=0$, and if $N_s<\infty$, then $\sigma_i=1/N_s > 0$.\\
\item[(v)] If $r(Q_u)\rightarrow 0$ as $\|u\|_{\E_1}\rightarrow \infty$, then $\sigma_s=\infty$.\\
\item[(vi)] If $r(Q_u)\le\xi$ for some $\xi>0$ and every solution $(n,u)\in\R^+\times\E_1^+$ to \eqref{30}, then $\sigma_i\ge 1/\xi$. In particular, if $r(Q_u)\le 1$ for every solution $(n,u)\in\R^+\times\E_1^+$ to \eqref{30}, then $\sigma_i=1$ corresponding to supercritical bifurcation.
\end{itemize}
\end{prop}

\begin{proof}
(i) We note that, according to \eqref{51}, any solution $(n,u)$ to \eqref{30} satisfies $u=\Pi_u(\cdot,0)u(0)$ with $u(0)=0$ if $n=0$, whence $u\equiv 0$ in this case. Hence $0\not\in \sigma$. 

(ii) Since $(1,0)\in \mathfrak{C}^+$ this is immediate.

(iii) Let $N_i=0$. Then there is a sequence $((n_j,u_j))$ in $\mathfrak{C}^+$ with $u_j\not\equiv 0$ and $r(Q_{u_j})\searrow 0$. From \eqref{72} we obtain $n_j\nearrow \infty$, whence $\sigma_s=\infty$. Conversely, let $\sigma_s=\infty$. Then there is a sequence $((n_j,u_j))$ in $\mathfrak{C}^+$ with $n_j\nearrow \infty$ and \eqref{72} ensures $r(Q_{u_j})\searrow 0$, whence $N_i=0$. Let now $N_i>0$. On the one hand, we necessarily have $1=n r(Q_u)\ge n N_i$ for all $(n,u)\in\mathfrak{C}^+$. Consequently, $n\le 1/N_i$ for $n\in\sigma$ so that $\sigma_s\le 1/N_i$. On the other hand, since $1=n r(Q_u)\le \sigma_s r(Q_u)$ for $(n,u)\in\mathfrak{C}^+$ we have $r(Q_u)\ge 1/\sigma_s$ for $u\in\Gamma$ and thus $N_i\ge 1/\sigma_s$.

(iv) The same as in (iii).

(v) Theorem \ref{A9} implies that $\sigma_s=\infty$ or there are $u_j\in\Gamma$ with $\|u_j\|_{\E_1}\rightarrow \infty$. Then $\sigma_s=\infty$ or $r(Q_u)\rightarrow 0$ by assumption, hence $\sigma_s=\infty$ in both cases according to (iii).

(vi) From (iv) we obtain $\sigma_i=1/N_s\ge 1/\xi$ if $N_s\le\xi$. In particular, if $\xi=1$, then $\sigma_i=1$ due to (ii).

\end{proof}

Of course, particularly interesting is the case $\sigma_s=\infty$. As in \cite[Cor.3]{Cushing2} one can easily impose conditions on $b$ and $\A=A+\mu$ guaranteeing $r(Q_u)\rightarrow 0$ as $\|u\|_{\E_1}\rightarrow 0$, whence $\sigma_s=\infty$ by (v) of the above proposition. In particular, as in \cite[Ex.3.3]{WalkerSIMA} one can put conditions on the same quantities leading to $r(Q_u)\le 1$, that is, to supercritical bifurcation in view of (vi).\\

As pointed out in Remark \ref{A8} the assumptions on the spectral radii and the properties of the eigenvalues of the operators $Q_u$ in Theorem \ref{A9} or Theorem \ref{A11} are not too restrictive but rather natural in applications due to the fact that these operators are compact and strongly positive in many cases.
In general, we refer to \cite[Thm.3.1, Ex.3.1, Ex.3.2, Ex.3.3]{WalkerSIMA} for other examples of diffusion operators $A$, birth moduli $b$, and death moduli $\mu$ satisfying the assumptions of Theorem~\ref{A9} or Theorem~\ref{A11} and provide here merely one example. Clearly, the set of applications is not restricted to following example.\\

\begin{example}\label{ex}

Let $\Om\subset\R^N$, $N\ge 1$, be a bounded and smooth domain lying locally on one side of $\partial\Om$. Let the boundary $\partial\Om$ be the distinct union of two sets $\Gamma_0$ and $\Gamma_1$ both of which are open and closed in $\partial\Omega$. For simplicity we assume $a_m\in (0,\infty)$ and consider a second order differential operator of the form
\bqn\label{50a}
\Ac(u,a,x)w:=-\nabla_x\cdot\big(D(a,x)\nabla_x w\big)+g\big(u(a),\nabla_xu(a)\big)\cdot\nabla_x w+h\big(u(a),\nabla_xu(a)\big)w\ ,
\eqn
where
\bqn
\begin{aligned}\label{50g}
&D:J\rightarrow C^1(\bar{\Om})\ \text{is bounded and uniformly H\"older}\\
&\text{ continuous with}\ D(a,x)>0\, ,\ (a,x)\in J\times \bar{\Omega}\ ,
\end{aligned}
\eqn
and
\bqn
\begin{aligned}\label{50aa}
&g\in C^3(\R\times\R^N,\R^N)\ \text{with}\ g(0,0)=0\ \text{and}\  h\in C^3(\R\times\R^N,\R)\ .
\end{aligned}
\eqn
Let 
\bqn\label{52}
\nu_0\in C^1(\Gamma_1)\ ,\qquad \nu_0(x)\ge 0\ ,\quad x\in\Gamma_1\ ,
\eqn
and let $\nu$ denote the outward unit normal to $\Gamma_1$. Let
\bqnn
\mathcal{B}(x)w:=\left\{\begin{array}{ll} w\ , & \text{on}\ \Gamma_0\ ,\\
 \frac{\partial}{\partial\nu}w+\nu_0(x) w\ , & \text{on}\ \Gamma_1\ .
\end{array}
\right.
\eqnn
Fix $p,q\in (1,\infty)$ with 
\bqn\label{50b}
\frac{1}{p}+\frac{N}{2q}<1\ ,
\eqn
and let $E_0:=L_q:=L_q(\Om)$ be ordered by its positive cone of functions that are nonnegative almost everywhere. Observe that 
$$
E_1:=\Wqb^2:=\Wqb^2(\Om):=\big\{u\in W_q^2\,;\, \mathcal{B}u=0\big\}\dhr L_q=E_0 
$$ 
and, up to equivalent norms, the interpolation spaces are subspaces of the Besov spaces $B_{q,p}^{2\xi}:=B_{q,p}^{2\xi}(\Om)$, that is,
\bqnn
E_\xi:=\big(L_q,W_{q,\mathcal{B}}^2\big)_{\xi,p}\,\dot{=}\,\Bqb^{2\xi}:=\left\{\begin{array}{ll} B_{q,p}^{2\xi}\ ,\quad &0<{2\xi}<1/q\ ,\\
\big\{w\in B_{q,p}^{2\xi}\,;\, u\vert_{\Gamma_0}=0\big\}\ ,& 1/q<{2\xi}<1+1/q\ , 2\xi\not= 1\ ,\\
\big\{w\in B_{q,p}^{2\xi}\,;\, \mathcal{B}u=0\big\}\ ,& 1+1/q<2\xi<2\ ,
\end{array}
\right.
\eqnn
(see e.g. \cite{Triebel}).
In particular, since $2(1-1/p)>N/q$, we have $E_{\varsigma}\,\dot{=}\,\Bqb^{2-2/p}\hookrightarrow C(\bar{\Om})$ for $\varsigma=1-1/p$
and thus $\mathrm{int}(E_\varsigma^+)\ne\emptyset$. Assume then further that
\bqn\label{50d}
\begin{aligned}
\mu:\R\times J\rightarrow [0,\infty)\ \text{is uniformly Lipschitz continuous on bounded sets}
\end{aligned}
\eqn
and if $\theta(a):=\mu(0,a)+h(0,0)$, then
\bqn\label{50e}
\begin{aligned}
& \theta : J\rightarrow (0,\infty)\ \text{is bounded and uniformly H\"older continuous}\ .
\end{aligned}
\eqn
Finally, let $b$ be such that
\bqn\label{50f}
b\in C^3(\R,(0,\infty))\ .
\eqn
Consider
$$
A(u,a)w:=\Ac(u,a,\cdot)w\ ,\quad w\in E_1\ ,\quad u\in \E_1\ ,
$$
and define $\A(u,a):=A(u,a) +\mu(u(a),a)$ for $u\in\E_1$ and $a\in J$. Then $$\A_0(a) w:=\A(0,a)w= -\nabla_x\cdot\big(D(a,\cdot)\nabla_x w\big)+\theta(a)w$$ is such that $\A_0:J\rightarrow \ml(\Wqb^2,L_q)$ is bounded and uniformly H\"older continuous by \eqref{50g} and \eqref{50e}. Moreover, from \cite[Sect.7,Thm.11.1]{AmannIsrael} it follows that for $a\in J$ fixed, $-\A_0(a)$ is resolvent positive, generates a contraction semigroup of negative type on each $L_r(\Om)$, $r\in (1,\infty)$, and is self-adjoint on $L_2(\Om)$. Hence $\A_0$ generates a positive parabolic evolution operator with regularity subspace $E_1=\Wqb^2$ by \cite[II.Cor.4.4.2]{LQPP} and possesses maximal $L_p$-regularity according to \cite[III.Ex.4.7.3,III.Thm.4.10.10]{LQPP}, whence \eqref{31}.
Owing to \eqref{50b} we may choose numbers $\bar{s}, s$, and $\alpha$ such that
\bqn\label{x}
\frac{1}{p}<\bar{s}<s<1-\alpha<\frac{1}{2}+\frac{1}{2p}-\frac{N}{4q}\ .
\eqn
Then, as in Remark \ref{KK}, we have
$$\E_1:=L_p(J,\Wqb^2)\cap W_p^1(J,L_q)\dhr W_p^s(J,\Bqb^{2\alpha})=:\F_j\ ,\quad j=1,\ldots,4\ ,$$ 
and analogously to \cite[Lem.2.7]{WalkerEJAM} we obtain that 
$$
[u\mapsto g(u,\nabla_xu)]: W_p^s(J,\Bqb^{2\alpha})\rightarrow W_p^{\bar{s}}(J,\Bqb^{2\alpha-1}) 
$$ 
is Lipschitz continuous. In particular, since $\bar{s}>1/p$ and $2\alpha-1>N/q$ we deduce for $u\in\E_1$ that $$[a\mapsto g(u(a),\nabla_x u(a))]:J\rightarrow C(\bar{\Om})$$ is H\"older continuous. Clearly, the same holds true for $h(u,\nabla_xu)$ and also
\bqn\label{xx}
[u\mapsto b(u)]: W_p^s(J,\Bqb^{2\alpha})\rightarrow W_p^{\bar{s}}(J,\Bqb^{2\alpha})\ \text{is Lipschitz continuous}\ . 
\eqn
Similarly we obtain from \eqref{50d} and the embedding $\E_1\hookrightarrow C^{\varsigma-\upsilon}(J,E_\upsilon)$, $\upsilon\le \varsigma$, (being due to the interpolation inequality \cite[I.Thm.2.11.1]{LQPP}) that $[a\mapsto \mu(u(a),a)] : J\rightarrow C(\bar{\Om})$ is H\"older continuous. Gathering these information and invoking \cite[II.Cor.4.4.2]{LQPP} we deduce that $\A(u,\cdot)$ generates for each $u\in\E_1$ a positive parabolic evolution operator $\Pi_u(a,\sigma)$, $0\le \sigma\le a<a_m$, on $E_0=L_q$, whence \eqref{51}. Obviously, \eqref{32} holds for $$\A_*(u,a):=\A(u,a)-\A_0(a)\ .$$ 
Note then that pointwise multiplication
$$
\Bqb^{2\alpha}\cdot \Bqb^{2\alpha}\hookrightarrow \Bqb^{2(1-1/p)}\doteq E_\varsigma
$$
is continuous according to \eqref{x} and \cite[Thm.4.1]{AmannMultiplication}. Thus, the same theorem together with $\bar{s}>1/p$ ensures that pointwise multiplication
$$
W_p^{\bar{s}}\big(J,\Bqb^{2\alpha}\big)\cdot W_p^{\bar{s}}\big(J,\Bqb^{2\alpha}\big)\hookrightarrow L_1\big(J,\Bqb^{2(1-1/p)}\big)
$$
is continuous since $J=(0,a_m)$ with $a_m<\infty$.
From this we obtain $\ell_0\in\ml(\E_1,E_1)$ and $\ell_*\in C(\F_4,E_\varsigma)$, where
$$
\ell_0(u):=b(0) \int_0^{a_m}u(a)\,\rd a\quad \text{and}\quad \ell_*(u):=\int_0^{a_m}[b(u(a))-b(0)]u(a)\,\rd a\ .
$$
Furthermore, 
$$
\|\ell_*(u)\|_{E_\varsigma}\le \| b(u)-b(0)\|_{W_p^{\bar{s}}(J,\Bqb^{2\alpha})} \|u\|_{W_p^{\bar{s}}(J,\Bqb^{2\alpha})}\ ,
$$ 
whence \eqref{34} and \eqref{35} since $s>\bar{s}$. Defining
$$
Q_u:=\int_0^{a_m}b(u(a))\Pi_u(a,0)\,\rd a\ ,\quad u\in\E_1\ ,
$$
it is immediate that $Q_u\in\K_+(E_\varsigma)$ due to \cite[II.Lem.5.1.3]{LQPP} and the compact embedding $E_1\dhr E_\varsigma$ since $a_m<\infty$. Moreover, $Q_u$ is strongly positive since $\Pi_u(a,0)$ is strongly positive on $E_\varsigma^+$ for $a\in J\setminus\{0\}$ (see \cite[Thm.13.6]{DanersKochMedina}). In particular, $Q_u$ is irreducible and so, by \cite[Thm.12.3]{DanersKochMedina}, $r(Q_u)>0$ is simple and the only eigenvalue of $Q_u$ with a positive eigenfunction since $\mathrm{int}(E_\varsigma^+)\ne\emptyset$ as observed above. This ensures \eqref{53} and \eqref{54}. Therefore, if $b$ is normalized such that $r(Q_0)=1$, we conclude thanks to Theorem~\ref{A9}:

\begin{prop}\label{vvvv}
Suppose \eqref{50a}-\eqref{50f} and let $r(Q_0)=1$. Then the problem
\begin{align*}
&\partial_a u+\Ac(u,a,x)u+\mu\big(u(a),a\big)u=0\ ,\quad a\in J\ ,\quad x\in\Om\ ,\\
&u(0,x)=n\int_0^\infty b\big(u(a,x)\big)u(a,x)\,\rd a\ ,\quad x\in\Om\ ,\\
&\mathcal{B}(x)u(a,x)=0\ ,\quad a>0\ ,\quad x\in\partial\Om\ ,
\end{align*}
admits an unbounded contiunuum $\mathfrak{C}^+$ of positive nontrivial solutions $(n,u)$ in  $$\R^+\times \big(L_p^+(J,\Wqb^2)\cap W_p^1(J,L_q)\big)\ .$$
\end{prop}

\end{example}

A noteworthy variant of the previous example is considering a functional (instead of a local) dependence of $A$, $b$, or $\mu$ on $u$ with respect to age. For details we refer to \cite[Ex.3.2, Ex.3.3]{WalkerSIMA}.

\section{Positive Solutions via a Fixed Point Argument}\label{Sec2}

The aim of this section is to give sufficient conditions for the existence of nontrivial nonnegative solutions to the (parameter-independent) problem
\bqn
\begin{aligned}\label{z}
&\partial_au \, +\, \A(u,a)\,u\, =0\ ,\qquad a\in J\ ,\\ 
&u(0)\, =\, \int_0^{a_m}b(u,a)\, u(a)\, \rd a\ ,
\end{aligned}
\eqn
in $E_0$ without assuming a decomposition \eqref{30A}-\eqref{32}. Due to the quasilinear structure of the first equation we require some assumptions that can considerably be weaken if one restricts to linear problems. For $\theta\in (0,1)$ we put $X_\theta:=L_1(J,E_\theta)$ and $X_\theta^+:=L_1^+(J,E_\theta)$. Let
\bqn\label{8}
0<\alpha<\beta<1\ .
\eqn
We suppose that, given any $R>0$, there are $\rho, \omega, \eta>0$, $\sigma\in\R$, and $\kappa\ge 1$ depending possibly on $R$ such that for  $\Phi_\alpha:=\mathbb{B}_{X_\alpha}(0,R)\cap X_\alpha^+$ we have
\bqn\label{9}
\begin{aligned} & \A\in C^{0,\rho}\big(\Phi_\alpha\times J,\ml(E_1,E_0)\big)\ \text{with}\\
&\big[\A(u,\cdot)\big]_{\rho,J}\le \eta\ ,\quad \sigma+\A(u,\cdot)\subset \mH(E_1,E_0;\kappa,\omega)\ ,\\
&\text{and}\ \A(u,\cdot) \ \text{is resolvent positive for each}\ u\in\Phi_\alpha\ .
\end{aligned}
\eqn
Observe that \eqref{9} and \cite[II.Cor.4.4]{LQPP} ensure that for each $u\in\Phi_\alpha$ there is a unique positive parabolic evolution operator $\Pi_u(a,\sigma)$, $0\le \sigma\le a<a_m$, on $E_0$ corresponding to $\A(u,\cdot)$. The evolution operator satisfies according to \cite[II.Lem.5.1.3]{LQPP} the estimates
\bqn\label{10}
\|\Pi_u(a,\sigma)\|_{\ml(E_\xi)}+(a-\sigma)^{\xi-\zeta/2}\|\Pi_u(a,\sigma)\|_{\ml (E_\zeta,E_\xi)}\le c_0e^{\nu(a-\sigma)}\ ,\quad 0\le\sigma\le a<a_m\ ,
\eqn
for $0\le\zeta\le\xi\le 1$ and some constants $c_0=c_0(R,\xi,\zeta)>0$ and $\nu=\nu(R)\in \R$ (independent of $u\in\Phi_\alpha$). We assume that
\bqn\label{11}
\nu<0\quad\text{if}\quad a_m=\infty\ .
\eqn
To control the dependence of the evolution operator $\Pi_u$ on $u\in\Phi_\alpha$ we require for each $u\in\Phi_\alpha$ the existence of $\varepsilon=\varepsilon(u)>0$ and a measurable function $g:(0,\varepsilon)\times J\rightarrow \R^+$ (depending possibly on $u$) with
\bqn\label{12}
\begin{aligned}
&\lim_{r\rightarrow 0^+}\int_0^{a_m}g(r,a)\,\rd a=0\ ,\\
& \max_{0\le\sigma\le a}\| \A(u,\sigma)-\A(\bar{u},\sigma)\|_{\ml(E_1,E_0)}\le g(\|u-\bar{u}\|_{X_\alpha^+},a)\ ,\quad a\in J\ , \quad \|u-\bar{u}\|_{X_\alpha^+}\le\varepsilon\ .
\end{aligned}
\eqn
As for $b$ appearing in \eqref{z} we suppose that 
\bqn\label{13}
\begin{aligned}
&b\in L_{\infty}^+(\Phi_\alpha\times J,F)\ ,\\
& e^{-\nu a}\| b(u,a)-b(\bar{u},a)\|_{F}\le g(\|u-\bar{u}\|_{X_\alpha^+},a)\ ,\quad a\in J\ , \quad \|u-\bar{u}\|_{X_\alpha^+}\le\varepsilon\ .
\end{aligned}
\eqn
Here $F$ is assumed to be a Banach space ordered by a convex cone $F^+$ such that a (bilinear) multiplication $m:=[(f,e)\mapsto fe]$ is induced which is continuous considered as mappings
\bqn\label{77}
m:F\times E_\beta\rightarrow E_\beta\quad\text{and}\quad m:F\times E_1\rightarrow E_\delta\quad\text{for some}\ \delta\in (\beta,1]
\eqn
and such that $m(f,e)=fe\in E_\beta^+$ for $f\in F^+$ and $e\in E_\beta^+$. Note that $F=\R$ is appropriate with $\delta=1$.

As a consequence of \eqref{10}, \eqref{11}, \eqref{13}, and the compact embedding $E_\delta\dhr E_\beta$ we have
\bqn\label{14}
Q(u)=\int_0^{a_m}b(u,a) \Pi_u(a,0)\,\rd a\in  \ml_+(E_\beta,E_\delta)\cap\K_+(E_\beta)\ ,\quad u\in\Phi_\alpha\ .
\eqn
Solutions to \eqref{z} are, as noted in the introduction, fixed points of the map $$(u,B)\mapsto \big(\Pi_u(\cdot,0)B,Q(u)B\big)$$
with $B=u(0)$.
Clearly, \eqref{8}-\eqref{77} are technical but not restrictive assumptions for applications (see \cite[Sect.3]{WalkerSIMA}). However, since the main task is to single out nontrivial solutions we also have to impose structural and thus more restrictive assumptions in order to apply a fixed point theorem for conical shells \cite{AmannSIAM}. The assumptions read:
\begin{align}\label{15}
&\text{there are}\ \tau_0>0\ \text{and}\ \psi\in E_\beta^+\ \text{with}\ \psi\not\in \bigcup_{\substack{u\in X_1^+\setminus\{0\}\\ \| u\|_{X_\alpha}<\tau_0}} \mathrm{rg}_+\big( 1-Q(u)\big)\ ,
\end{align}
where $\mathrm{rg}_+( 1-Q(u)):=\{(1-Q(u))B\,;\, B\in E_\beta^+\}$, and
\begin{align}\label{16}
&\text{there is}\ \tau_1>0\ \text{such that}\ r(Q(u))\le 1\ \text{for}\ u\in X_1^+\ \text{with}\ \|u\|_{X_\beta}\ge \tau_1\ ,
\end{align}
where $r(Q(u))$ denotes the spectral radius of the operator $Q(u)\in\ml(E_\beta)$.\\

We comment in more detail on the structural requirements \eqref{15}, \eqref{16} after the proof of the following result, which is in the spirit of \cite[Thm.1]{Pruess2}:

\begin{thm}\label{A1}
Suppose \eqref{8}, \eqref{9}, \eqref{11}-\eqref{77}, \eqref{15}, and \eqref{16}. Then \eqref{z} has at least one nontrivial nonnegative solution 
$$u\in L_1(J,E_1)\cap C^1(J\setminus\{0\},E_0)\cap C(J, E_\delta)\ . $$
\end{thm}

\begin{proof}
We shall employ \cite[Thm.12.3]{AmannSIAM} in proving the statement. Let $\X:=X_\alpha^+ \times E_\beta^+$ and $\X_R:=\mathbb{B}_{\X}(0,R)$ with $$R:=\tau_1\big(\|i\|_{\ml(E_\beta,E_\alpha)}+\|b\|_{L_\infty(\mathbb{B}_{X_\beta^+}(0,\tau_1)\times J,F)}\big)>0\ ,
$$
where $i$ is the natural injection $E_\beta\hookrightarrow E_\alpha$.
We put
$$
f(u,B):=\big(\Pi_u(\cdot,0)B,Q(u)B\big)
$$ and first claim that $f:\X_R\rightarrow \X$ is continuous and $f(\X_R)$ is relatively compact in $\X$. Indeed, given $(u,B), (\bar{u},\bar{B})\in\X_R$ we note that
\bqnn
\begin{split}
\int_0^{a_m}\| \Pi_u(a,0)B-\Pi_{\bar{u}}(a,0)\bar{B}\|_{E_\alpha}\, \rd a &\le  \int_0^{a_m}\|\Pi_u(a,0)-\Pi_{\bar{u}}(a,0)\|_{\ml(E_\alpha)}\, \rd a\, \|B\|_{E_\alpha} \\ 
&\qquad +\int_0^{a_m}\|\Pi_{\bar{u}}(a,0)\|_{\ml(E_\alpha)}\, \rd a\, \|B-\bar{B}\|_{E_\alpha}
\\
&\le c(R)\int_0^{a_m}g\big(\|u-\bar{u}\|_{X_\alpha^+},a\big)\, \rd a\, +\, c(R)\| B-\bar{B}\|_{E_\alpha}\ ,
\end{split}
\eqnn
where we invoked \eqref{10}-\eqref{12}, and \cite[II.Lem.5.1.4]{LQPP}. Thus $\Pi_{\bar{u}}(a,0)\bar{B}\rightarrow \Pi_u(a,0)B$ in $X_\alpha$ as $(\bar{u},\bar{B})$ approaches $(u,B)$ in $\X_R$ by \eqref{12}. Similarly we deduce $Q(\bar{u})\bar{B}\rightarrow Q(u)B$ in $E_\beta$ as $(\bar{u},\bar{B})\rightarrow (u,B)$ in $\X_R$, whence the continuity of $f$. Next, we use the characterization for compact sets in $X_\alpha=L_1(J,E_\alpha)$ due to \cite{Gutman}. We may assume $a_m=\infty$. The previous argument entails
$$
\sup_{(u,B)\in\X_R}\| \Pi_u(\cdot,0)B\|_{X_\alpha}<\infty\ .
$$
Moreover,
$$
\int_N^\infty  \|\Pi_u(a,0)B\|_{E_\alpha}\, \rd a\, \le\, c(R)\int_N^\infty e^{\nu a}\, \rd a\,\longrightarrow 0\quad\text{as}\quad N\longrightarrow \infty
$$
uniformly with respect to $(u,B)\in \X_R$ by \eqref{10} and \eqref{11}. Let $h>0$. Then, from \eqref{10}, \eqref{11}, and equation (II.5.3.8) in \cite{LQPP} we deduce
\bqnn
\begin{split}
\int_0^{\infty}\| \Pi_u(a+h,0)B-\Pi_{u}(a,0)B\|_{E_\alpha}\, \rd a\, &\le\,  \int_0^{\infty}\| \Pi_u(a+h,0)-\Pi_{u}(a,0)\|_{\ml(E_\beta,E_\alpha)}\, \rd a\, \|B\|_{E_\beta}\\
& \le\, c(R)\, h^{\beta-\alpha}
\end{split}
\eqnn
and the right hand side tends to 0 as $h\rightarrow 0$ uniformly with respect to $(u,B)\in \X_R$ in view of \eqref{8}. Furthermore, since \eqref{10}, \eqref{11} ensure
$$\|\Pi_u(a,0)B\|_{E_\beta}\le c(R)\ , \quad a\in (0,\infty)\ ,\quad (u,B)\in \X_R\ ,
$$
we obtain from the compact embedding $E_\beta\dhr E_\alpha$ that $\Pi_u(a,0)B$ belongs to a fixed compact subset of $E_\alpha$. Applying now \cite[Thm.A.1]{Gutman} we derive the relative compactness of the set $\{\Pi_u(\cdot,0)B\,;\, (u,B)\in \X_R\}$ in $X_\alpha$. Next, observing that
$$
\|Q(u)B\|_{E_\delta}\le\int_0^{a_m} \| b(u,a)\|_F\, \|\Pi_u(a,0)\|_{\ml(E_\beta,E_1)}\,\rd a\,\| B\|_{E_\beta}\le c(R)
$$
for $(u,B)\in\X_R$ according to \eqref{10}, \eqref{13}, and \eqref{77} we may use the compact embedding $E_\delta\dhr E_\beta$ and also obtain the relative compactness of the set $\{Q(u)B\,;\, (u,B)\in \X_R\}$ in $E_\beta$. Therefore, $f(\X_R)$ is relatively compact in $\X$. It remains to check the crucial conditions (i) and (ii) from \cite[Thm.12.3]{AmannSIAM}. For (i) suppose there exist $\lambda>1$ and $(u,B)\in \X_R$ for which
$$
\|u\|_{X_\alpha}+\|B\|_{E_\beta}=\|(u,B)\|_{\X_R}=R\quad\text{and}\quad f(u,B)=\lambda (u,B)\ ,
$$
that is,
\bqn
\begin{aligned}
\lambda u(a)&=\Pi_u(a,0)B\ ,\quad a\in J\ ,\label{22}\\
\lambda B&=Q(u)B\ . 
\end{aligned}
\eqn
Since $\lambda>1$ we have $B\not=0$ (otherwise $u\equiv 0$ contradicting $R>0$). From \eqref{8}, \eqref{10}, \eqref{11}, \eqref{14} we deduce, on the one hand, that $B\in E_\delta^+$ is an eigenvector for $Q(u)$ corresponding to the eigenvalue $\lambda>1$ and, on the other hand, that $u\in X_1^+$. Invoking \eqref{16} we see that this is only possible if
$\|u\|_{X_\beta}<\tau_1$. 
Consequently, recalling \eqref{77} and \eqref{22} we derive the contradiction
\bqnn
\begin{split}
R&=\|u\|_{X_\alpha}+\|B\|_{E_\beta}=\int_0^{a_m}\| u(a)\|_{E_\alpha}\,\rd a\,+\,\frac{1}{\lambda}\left\|\int_0^{a_m} b(u,a)\Pi_u(a,0)B\,\rd a\right\|_{E_\beta}\\
&<\tau_1\|i\|_{\ml(E_\beta,E_\alpha)}+\left\|\int_0^{a_m} b(u,a)u(a)\,\rd a\right\|_{E_\beta}\\
&\le\tau_1\|i\|_{\ml(E_\beta,E_\alpha)}\,+\,\|b\|_{L_\infty(\mathbb{B}_{X_\beta^+}(0,\tau_1)\times J,F)}\,\|u\|_{X_\beta}\,<\, R\ .
\end{split}
\eqnn
This  ensures $f(u,B)\not= \lambda(u,B)$ for all $\lambda>1$ and all $(u,B)\in \X_R$ with $\|(u,B)\|_{\X_R}=R$, whence (i) from \cite[Thm.12.3]{AmannSIAM}. Finally, let $\psi$ be as in \eqref{15} with $\tau_0<R$ and assume there exists $\lambda>0$ and $(u,B)\in \X_R$ with $\|(u,B)\|_{\X_R}=\tau_0$ and $(u,B)-f(u,B)=\lambda (0,\psi)$. Then $u=\Pi_u(\cdot,0)B$, hence $u\in X_1$ by \eqref{10}, \eqref{11} with $\|u\|_{X_\alpha}<\tau_0$, and $B=Q(u)B+\lambda\psi$. The latter implies $\psi\in\mathrm{rg}_+(1-Q(u))$ contradicting \eqref{15}. Thus (ii) from \cite[Thm.12.3]{AmannSIAM} is verified, too, and we conclude a fixed point $(u,B)\in \X_R\setminus\{(0,0)\}$ of the map $f$, that is, a nontrivial positive solution to \eqref{z}. As for the additional regularity stated in the theorem we observe that necessarily $u(a)=\Pi_u(a,0)B$ for $a\in J$ with $B=Q(u)B\in E_\delta$. It thus suffices to refer to the regularity theory of Chapter II in \cite{LQPP}. 
\end{proof}

Example \ref{ex} or the examples in \cite[Sect.3]{WalkerSIMA} apply with minor modifications to the situation of Theorem~\ref{A1}. We note that the special assumptions \eqref{11}, \eqref{15}, and \eqref{16} are also not too hard to verify in applications in view of the following remark.

\begin{rems}\label{A2}
(a) Let $\A(u,a)$ be of the form \eqref{a} with $\mu$ being real-valued so that its evolution operator is given by $$\Pi_u(a,\sigma)=e^{-\int_\sigma^a \mu(u,r)\rd r}U_{A(u,\cdot)}(a,\sigma)\ , $$ where $U_{A(u,\cdot)}$ denotes the parabolic evolution operator corresponding to $A(u,\cdot)$. Then \eqref{11} holds provided that there is $\mu_0\ge 0$ such that $\varliminf\limits_{a\rightarrow \infty}\mu(u,a)\ge \mu_0$ uniformly with respect to $u\in\Phi_\alpha$ and $\mu_0>s(-A(u,a))$ for $u\in\Phi_\alpha$ and $a\in (0,\infty)$ with $s(-A(u,a))$ being the spectral bound of the operator $-A(u,a)$ considered as a linear operator in $E_0$ (see \cite[Sect.I.1, Sect.II.5]{LQPP}).\\

(b) Suppose  \eqref{14}. Then condition \eqref{16} is equivalent to assume that $\mathrm{ker}\big(\lambda-Q(u)\big)\cap E_\beta^+=\{0\}$ for all $\lambda>1$ and $u\in X_1^+$ with $\|u\|_{X_\beta}\ge \tau_1$.

\begin{proof} This follows from \eqref{14} and the Krein-Rutman theorem which states that $r(Q(u))>0$ is an eigenvalue of $Q(u)\in\K_+(E_\beta)$ with a positive eigenvector.
\end{proof}

(c) Suppose  \eqref{14}. Then condition \eqref{16} holds if $\|Q(u)\|_{\ml(E_\theta)}\le 1$ for some $\theta\in [0,\delta]$ and all $u\in X_1^+$ with $\|u\|_{X_\beta}\ge \tau_1$.

\begin{proof} 
This is a consequence of (b) and \eqref{14}. 
\end{proof}

(d) Suppose  \eqref{14}. Then condition \eqref{15} is satisfied provided
$$ Q(u)-1\in\ml_+(E_\beta)\quad\text{for}\ u\in X_1^+\setminus\{0\}\ \text{with}\ \| u\|_{X_\alpha}<\tau_0\ .
$$
Note that the latter condition corresponds in the non-diffusive case $A\equiv 0$ to assuming the scalar inequality $Q(u)\ge 1$ for $\vert u\vert$ small as in \cite[Thm.1]{Pruess2} and \cite[Thm.4.1]{Webb}.

\begin{proof} 
Since in this case $\mathrm{rg}_+(1-Q(u))\subset -E_\beta^+$ we may choose $\psi\in E_\beta^+\setminus\{0\}$ arbitrarily. 
\end{proof}

(e) Suppose  \eqref{14}. Then condition \eqref{15} holds provided there is $\psi\in \mathrm{ker}(1-Q(u))\cap E_\beta^+\setminus\{0\}$ such that $Q(u)\in \K_+(E_\beta)$ is irreducible for each $u\in X_1^+\setminus\{0\}$ with $\| u\|_{X_\alpha}<\tau_0$ and, e.g., the interior of $E_\beta^+$ is nonempty. 

\begin{proof} If $\psi$ is as in the statement, then the Krein-Rutman theorem (e.g. see \cite[Thm.12.3]{DanersKochMedina}) warrants that $r(Q(u))=1$ is a simple eigenvalue of $Q(u)$, hence $$\mathrm{ker}(1-Q(u))\cap\mathrm{rg}(1-Q(u))=\{0\}\ ,$$ from which we conclude $\psi\not\in\mathrm{rg}_+(1-Q(u))$ for each $u$. 
\end{proof}

(f) If effects of small populations are negligible, then \eqref{15} holds. More precisely, let $Q(u)=Q(0)$ for small $u\in X_1^+$ and let $Q(0)\in \K_+(E_\beta)$ be irreducible, $r(Q(0))=1$, and $\mathrm{int}(E_\beta^+)\not=\emptyset$. Then there is $\psi\in E_\beta^+\setminus\{0\}$ with $\psi\in\mathrm{ker}(1-Q(0))$ according to the Krein-Rutman theorem \cite[Thm.12.3]{DanersKochMedina}, and \eqref{15} follows from (e).\\

(g) Suppose  \eqref{14}. Then condition \eqref{15} is satisfied if there are $\tau_0>0$ and $\psi\in\mathrm{int}( E_\beta^+)$ such that, given any $u\in X_1^+\setminus\{0\}$ with $\| u\|_{X_\alpha}<\tau_0$, $Q(u)$ is irreducible and there exists $\alpha(u,\cdot)\in L_{1,loc}(0,a_m)$ with
$$
-\A(u,a)\psi=\alpha(u,a)\psi\ ,\quad a\in J\ , \qquad\text{and}\qquad \int_0^{a_m}b(u,a)e^{\int_0^a \alpha(u,r)\rd r}\,\rd a=1\ .
$$
It is worthwhile to remark, however, that in this case a nontrivial solution to \eqref{z} can be found also in the form $u(a)=\phi(a)\psi$, where the existence of a nonnegative nontrivial $\phi$ follows from \cite{Pruess2}.

\begin{proof} The assumptions imply $\Pi_u(a,0)\psi=e^{\int_0^a \alpha(u,r)\rd r}\psi$, $a\in J$, for $u\in X_1^+\setminus\{0\}$ with $\| u\|_{X_\alpha}<\tau_0$, whence $\psi\in\mathrm{ker}(1-Q(u))$, and we may apply (e).
\end{proof}

(h) Clearly, Theorem \ref{A1} applies to models involving several species, say with densities $u_j$, $1\le j\le N$, and $u=(u_1,\ldots,u_N)$. If $\A$ and $b$ in \eqref{z} have ``diagonal form'', that is, if each $u_j$ satisfies \eqref{z} with $\A(u,a)$ and $b(u,a)$ replaced by $\A_j(u,a)$ and $b_j(u,a)$, respectively, then it suffices to assume \eqref{15} and \eqref{16} for some component. More precisely, the assertion of Theorem~\ref{A1} holds true provided there are $j\in\{1,\ldots,N\}$ (at least one) and \eqref{15} is replaced by
\bqnn
\begin{aligned}
&\text{there are}\ \tau_0>0,\ \text{and}\ \psi_j\in E_\beta^+
 \text{with}\ \psi_j\not\in \bigcup_{\substack{u\in (X_1^+)^N\\ \| u\|_{(X_\alpha)^N}<\tau_0}} \mathrm{rg}_+\big( 1-Q_j(u)\big)\ 
\end{aligned}
\eqnn
while \eqref{16} is replaced by
\begin{align*}
&\text{there is}\ \tau_1>0\ \text{such that}\ r(Q_j(u))\le 1\ \text{for}\ u\in (X_1^+)^N\ \text{with}\ \|u\|_{(X_\beta)^N}\ge \tau_1\ ,
\end{align*}
where
$$Q_j(u):=\int_0^{a_m}b_j(u,a)\Pi_{u,j}(a,0)\,\rd a
$$
with $\Pi_{u,j}$ denoting the parabolic evolution operator corresponding to $\A_j(u,\cdot)$.

\begin{proof}
Looking for solutions $u=(u_1,\ldots,u_N)$, where only the $j$-components are non-vanishing, this follows by an obvious modification of the proof of Theorem \ref{A1}.
\end{proof}

\end{rems}

\end{document}